%% file: DeWit99d.tex
\documentclass[10pt]{article}

\usepackage{amssymb}

\title{
  Four Easy Pieces -- Explicit {R} Matrices from the
  $(\dot{0}_m|\alpha)$ Highest Weight Representations of
  {$U_q[gl(m|1)]$}
}

\author{David {De Wit}%
  \footnote{%
    RIMS, Kyoto University, 606-8502 Japan.
    \texttt{ddw@kurims.kyoto-u.ac.jp}
  }
}

\begin{document}

\maketitle

\begin{abstract}
   We provide explicit presentations of members of a suite of R
   matrices arising from the $(\dot{0}_m|\alpha)$ representations of
   the quantum superalgebras $U_q[gl(m|1)]$.  Our algorithm constructs
   both trigonometric and quantum R matrices; all of which are
   \emph{graded}, in that they solve a graded Yang--Baxter equation.
   This grading is easily removed, yielding R matrices that solve the
   usual Yang--Baxter equation.  For $m>2$, the computations are
   impracticable for a human to perform, so we have implemented the
   entire process in \textsc{Mathematica}, and then performed the
   computations for $m=1,2,3$ and $4$.
\end{abstract}


\section{Overview}

\noindent
This paper describes the results of the automation of an algorithm to
explicitly generate several R matrices.  Specifically, we construct
trigonometric R matrices $\check{R}^m (u)$ corresponding to
the $\alpha$-parametric highest weight minimal representations labeled
$(\dot{0}_m|\alpha)$, of the quantum superalgebras $U_q[gl(m|1)]$.
These representations are $2^m$ dimensional, irreducible,
and contain free complex parameters $q$ and $\alpha$; the real variable
$u$ is a spectral parameter.  Quantum R matrices $\check{R}^m$ are
immediately obtainable as the spectral limits $u\to\infty$ of the
$\check{R}^m (u)$.

Our R matrices are in fact \emph{graded}, as they are based on graded
vector spaces, hence they actually satisfy graded Yang--Baxter
equations.  However, it is a simple matter to remove this grading and
transform them into objects that satisfy the usual Yang--Baxter
equations.

\enlargethispage{\baselineskip}

These R matrices are of physical interest in that they are applicable
to the construction of exactly solvable models of interacting
fermions.  Corresponding to $\check{R}^m (u)$, we may construct an
integrable $2^m$ state fermionic model on a lattice.  Models associated
with $m=2$ and $m=3$ have been discussed in
\cite{GouldHibberdLinksZhang:96} and \cite{GeGouldZhangZhou:98a},
respectively. The $m=4$ case has an elegant interpretation in terms of
a $2$-leg ladder model for interacting electrons: a discussion of
this is provided in \S\ref{sec:Application}.

Furthermore, from each of our $\check{R}^m$, we may obtain a two
variable polynomial `Links--Gould' link invariant $LG^m$
\cite{LinksGould:92b}.  $LG^1$ degenerates to being the
Alexander--Conway polynomial in the single variable $q^{2\alpha}$
(c.f.  \cite{BrackenGouldZhangDelius:94b}).  $LG^2$ is in fact more
powerful than the well known two variable HOMFLY and Kauffman
invariants, although it cannot distinguish mutants or inversion
\cite{DeWit:99a,DeWitKauffmanLinks:99a}.  $LG^m$ for $m>2$ have similar
gross properties to $LG^2$, although they should be able to distinguish
more links \cite{DeWit:99e}.

The process has been implemented in \textsc{Mathematica}, and R
matrices computed for $m=1,2,3$ and $4$.  The description of the
computational details of the algorithms used to construct the R
matrices is rather long, and will be provided elsewhere
\cite{DeWit:99c}, as will observations relating to the construction of
the new invariants \cite{DeWit:99e}.


\section{Algebraic Details}

Fixing $m$, we are initially interested in a $2^m$ dimensional vector
space $V$ that is a module for the $U_q[gl(m|1)]$ representation
$\Lambda=(\dot{0}_m|\alpha)$.  The algebra contains a free complex
variable $q$, whilst the representation $\pi_{\Lambda}$ acting on $V$
contains a free complex variable $\alpha$.  Our $V$ is actually
($\mathbb{Z}_2$) \emph{graded}; this ensures compatibility with the
($\mathbb{Z}_2$) grading of $U_q[gl(m|1)]$. A full description of
$U_q[gl(m|n)]$ in terms of generators and relations is contained in
\cite[pp1237-1238]{Zhang:93}; for our purposes a set of \emph{simple}
generators for $U_q[gl(m|1)]$ is:
\begin{eqnarray*}
  \left\{
    \begin{array}{lll}
      K_a,           & a = 1, \dots, m+1 & \mathrm{(Cartan)},
      \\
      {E^{a}}_{a+1}, & a = 1, \dots, m   & \mathrm{(raising)},
      \\
      {E^{a+1}}_{a}, & a = 1, \dots, m   & \mathrm{(lowering)}
    \end{array}
  \right\}.
\end{eqnarray*}

We apply the Kac induced module construction (KIMC) \cite{Kac:78} to
establish a \emph{weight} basis $\{v_i\}_{i=1}^{2^{m}}$ for $V$.  This
involves postulating $v_1$ as a highest weight vector, and recursively
acting on $v_1$ with all possible distinct products of simple lowering
generators ${E^{a+1}}_a$ to define the other basis vectors, normalising
as we go.  This construction requires a PBW basis for $U_q[gl(m|1)]$,
which enables us to transform any element of the algebra into a normal
form (\cite{Zhang:93}, see also \cite{DeWit:99c}).

The tensor product module $V\otimes V$ has a natural (weight) basis
$\{v_i \otimes v_j\}_{i,j=1}^{2^m}$, which inherits the grading of
$V$.  To build R matrices acting on $V\otimes V$, we require an
alternative, orthonormal weight basis $B$ for $V \otimes V$
corresponding to its decomposition into irreducible $U_q[gl(m|1)]$
submodules. Again using the KIMC, the basis vectors of $B$ are derived
as linear combinations of the form $\gamma_{ij}(v_i\otimes v_j)$, where
the coefficients $\gamma_{ij}$ are algebraic expressions in $q$ and
$\alpha$.  (This process initially yields a basis for each submodule
that is not necessarily orthonormal, so we also apply a Gram--Schmidt
process.)

\pagebreak

For our particular representation, the orthogonal decomposition of
$V\otimes V$ contains no multiplicities
\cite[(34)]{DeliusGouldLinksZhang:95b}:
\begin{eqnarray*}
  V \otimes V
  & \triangleq &
  \bigoplus_{k=1}^{m+1}
    V_k,
\end{eqnarray*}
where $V_k$ has highest weight
$\lambda_k=(\dot{0}_{m+1-k}, \dot{-1}_{k-1}\,|\,2\alpha+k-1)$.

The R matrices are then formed as weighted sums of projectors onto
these submodules $V_k$.  Explicitly, where $\check{P}_k$
is the projector onto submodule $V_k$, we have:
\begin{eqnarray*}
  \check{R}^m (u)
  = 
  \sum_{k=1}^{m+1}
    \Xi_k
    \check{P}_k,
  \qquad \qquad
  \check{R}^m
  =
  \sum_{k=1}^{m+1}
    \xi_k
    \check{P}_k,
\end{eqnarray*}
where $\Xi_k$ and $\xi_k$ are the following eigenvalues
of the R matrices on the submodules $V_k$
(\cite{GouldLinksZhang:96b}, and c.f. \cite{GeGouldZhangZhou:98a}):
\begin{eqnarray*}
  \Xi_k
  & = &
  \prod_{j=0}^{k-2}
    \frac{
      [\alpha + j + u]_q
    }{
      [\alpha + j - u]_q
    },
  \\
  \xi_k
  & = &
  \lim_{u \to \infty}
    \Xi_k
  =
  {(-1)}^{k-1}
  q^{(k - 1) (2 \alpha + k - 2)},
\end{eqnarray*}
where we intend $\Xi_1=\xi_1=1$, and we have used
the $q$ bracket:
\begin{eqnarray*}
  {[X]}_q
  \triangleq
  \frac{q^{+ X} - q^{- X}}{q^{+ 1} - q^{- 1}};
  \qquad
  \mathrm{observe~that}
  \qquad
  \lim_{q \to 1}{[X]}_q
  =
  X.
\end{eqnarray*}
Thus we have the intended spectral limit
$\check{R}^m = \lim_{u\to\infty}\check{R}^m (u)$.
The resulting R matrices are normalised such that the coefficients of
the `first' components (viz $e^{1 1}_{1 1}$) are unity.  For
applications, other choices of normalisation may be applicable
\cite{DeWit:99e}.

To be certain, $\check{R}^m (u)$ satisfies the following
graded version of the trigonometric Yang--Baxter equation:
\begin{equation}
   \begin{array}{l}
     {(-1)}^{[b'][c'] + [a'][c] + [a][b] + [b'][b'']}
     {\check{R} (u)    }^{c'' b''}_{b' c'}
     {\check{R} (u + v)}^{c'  a''}_{a' c }
     {\check{R} (v)    }^{b'  a' }_{a  b }
     \\[1mm]
     \qquad
     \qquad
     =
     {(-1)}^{[a'][b'] + [a][c'] + [b][c] + [b][b']}
     {\check{R} (v)    }^{b'' a''}_{a' b'}
     {\check{R} (u + v)}^{c'' a' }_{a  c'}
     {\check{R} (u)    }^{c'  b' }_{b  c },
   \end{array}
   \label{eq:ComponentGTYBE}
\end{equation}
where $[a]$ is the grading of the vector $v_a$.
The parity factors in
(\ref{eq:ComponentGTYBE}) may be removed by the following
transformation (e.g. see \cite{DeWit:98}):
\begin{eqnarray*}
  \begin{array}{ll}
    \check{R}^{a' b'}_{a b} (u)
    \mapsto
    {(-1)}^{[a]([b]+[b'])}
    \check{R}^{a' b'}_{a b} (u),
  \end{array}
\end{eqnarray*}
after which $\check{R}(u)$ which satisfied (\ref{eq:ComponentGTYBE})
now satisfies the usual ungraded TYBE:
\begin{equation}
  {\check{R} (u)    }^{c'' b''}_{b' c'}
  {\check{R} (u + v)}^{c'  a''}_{a' c }
  {\check{R} (v)    }^{b'  a' }_{a  b }
  =
  {\check{R} (v)    }^{b'' a''}_{a' b'}
  {\check{R} (u + v)}^{c'' a' }_{a  c'}
  {\check{R} (u)    }^{c'  b' }_{b  c },
  \label{eq:ComponentTYBE}
\end{equation}
written in noncomponent form as:
\begin{equation}
  \check{R}_{12} (u)
  \check{R}_{23} (u + v)
  \check{R}_{12} (v)
  =
  \check{R}_{23} (v)
  \check{R}_{12} (u + v)
  \check{R}_{23} (u).
  \label{eq:NoncomponentTYBE}
\end{equation}

\pagebreak

\noindent
In the spectral limit
$\check{R}=\lim_{u\to\infty}\check{R} (u)$,
this of course becomes a QYBE:
\begin{equation}
  \check{R}_{12}
  \check{R}_{23}
  \check{R}_{12}
  =
  \check{R}_{23}
  \check{R}_{12}
  \check{R}_{23},
  \label{eq:NoncomponentQYBE}
\end{equation}
viz
$
  (\check{R} \otimes I)
  (I \otimes \check{R})
  (\check{R} \otimes I)
  =
  (I \otimes \check{R})
  (\check{R} \otimes I)
  (I \otimes \check{R})
$,
familiar as the braid relation
$\sigma_1 \sigma_2 \sigma_1 = \sigma_2 \sigma_1 \sigma_2$.

Defining $R(u) \triangleq P \check{R}(u)$, where
$P$ is a permutation operator, yields a trigonometric R matrix
$R(u)$ satisfying the following version of
(\ref{eq:NoncomponentTYBE}):
\begin{equation}
  R_{12} (u)
  R_{13} (u+v)
  R_{23} (v)
  =
  R_{23} (v)
  R_{13} (u+v)
  R_{12} (u).
  \label{eq:NoncomponentAlternativeTYBE}
\end{equation}
This transformation amounts to the mapping:
$
  {R (u)}^{a' b'}_{a b}
  =
  {\check{R} (u)}^{b' a'}_{a b}
$.
In component form, (\ref{eq:NoncomponentAlternativeTYBE}) is
more symmetric than (\ref{eq:ComponentTYBE}):
\begin{eqnarray*}
  {R (u)  }^{b'' c''}_{b' c'} 
  {R (u+v)}^{a'' c' }_{a' c } 
  {R (v)  }^{a'  b' }_{a  b }
  =
  {R (v)  }^{a'' b''}_{a' b'} 
  {R (u+v)}^{a'  c''}_{a  c'} 
  {R (u)  }^{b'  c' }_{b  c }.
\end{eqnarray*}


\section{Implementation}

The entire process has been implemented as a suite of functions in the
interpreted environment of \textsc{Mathematica}.  Whilst there is no
theoretical limit to $m$, storage and patience mean that a current
reasonable practical limit for $m$ is $4$.  The computations are
computationally inefficient!  Translation of the several thousand lines
of \textsc{Mathematica} code into a compiled language would increase
the speed of the algorithm enormously, but storage requirements would
still limit $m$.


\section{Results}

Both $\check{R}^m (u)$ and $\check{R}^m$ have been obtained
for $m = 1,2,3,4$.  Of these, the $m=1$ case (c.f.
\cite{BrackenGouldZhangDelius:94b}), can be done by hand in a couple of
hours; the complete $m=2$ case appears in my PhD thesis
\cite{DeWit:98}, and took several weeks to do by hand; partial details
of the $m=3$ case appear in \cite{GeGouldZhangZhou:98a}; whilst the
$m=4$ case is new.
By direct substitution, we have been able to verify that each
$\check{R}^m (u)$ satisifes (\ref{eq:NoncomponentTYBE})%
\footnote{%
  We didn't check $\check{R}^4 (u)$, as the computations
  would have been excessively expensive.
}
and that each $\check{R}^m$ satisfies
(\ref{eq:NoncomponentQYBE}).

Fixing $m$, each R matrix contains $2^{4m}$ (albeit mostly zero)
components.  Let $N'_m$ and $N_m$ be the number of nonzero
components of $\check{R}^m(u)$ and $\check{R}^m$
respectively.  As $\check{R}^m$ is the spectral limit of
$\check{R}^m (u)$, we expect $N_m \leqslant N'_m$.
The numbers of nonzero components of each of the $m+1$ projectors are
similar to $N_m $ and $N'_m$; and $N'_m = 6^m$
(why?). Let $s_m \triangleq N_m/2^{4m}$ be the sparsity of
$\check{R}^m$.  Table \ref{tab:data} summarises our results.

\begin{table}[ht]
  \begin{centering}
  \begin{tabular}{crcrrr}
    $m$ &$2^{4m}$ &       projector sizes    & $N'_m$ & $N_m$ & $s_m (\%)$ \\[1mm]
    \hline
        &         &                          & \\[-3mm]
    $1$ & $   16$ & $          5,5         $ & $   6$ & $  5$ & $31.3$ \\
    $2$ & $  256$ & $       25,34,25       $ & $  36$ & $ 26$ & $10.2$ \\
    $3$ & $ 4096$ & $   125,199,199,125    $ & $ 216$ & $139$ & $ 3.4$ \\
    $4$ & $65536$ & $625,1124,1254,1124,625$ & $1296$ & $758$ & $ 1.1$ \\
  \end{tabular}
  \caption{Numbers of nonzero components of our R matrices.}
  \label{tab:data}
  \end{centering}
\end{table}

Listings of the nonzero components of our R matrices are supplied in
the Appendix.

\vfill

\pagebreak


\section{An Application}
\label{sec:Application}

Of particular new interest is the interpretation of our $U_q[gl(4|1)]$
trigonometric R matrix $\check{R}^4(u)$ in the construction of an
exactly solvable $2$-leg ladder model of interacting electrons.
To this end, consider a $2$-leg ladder, with electron
occupation sites at the end of each rung.  Each site may contain a
maximum of $2$ electrons, each in state spin up $\uparrow$ or down
$\downarrow$. Thus, at each site we have $4$ possible states:
unoccupied $\vert 0\rangle$, both up
$\vert\!\!\uparrow\uparrow\rangle$, both down
$\vert\!\!\downarrow\downarrow\rangle$, or mixed
$\vert\!\!\uparrow\downarrow\rangle$. Taken together, each rung space,
corresponding to our $V$, is $16$ dimensional. These $16$ dimensions
correspond to $1$ electron-free state, $4$ single-electron states, $6$
two-electron states, $4$ three-electron states and $1$ four-electron
state.

A Hamiltonian determined by $\check{R}^4(u)$ describes the interactions
between rungs.  Discernment of the details of the terms contained
within this Hamiltonian are left as an exercise for the reader with some
idle time; the procedure essentially follows
\cite{%
  BrackenGouldZhangDelius:94b,%
  GeGouldZhangZhou:98a,%
  GouldHibberdLinksZhang:96%
}.


\section*{Acknowledgements}

My research at Kyoto University is funded by a Postdoctoral Fellowship
for Foreign Researchers (\# P99703), provided by the Japan Society for
the Promotion of Science.  D\={o}mo arigat\={o} gozaimashita!

Some of this work was completed under the direction of Dr Mark Gould as
part of ongoing research at The University of Queensland, Australia.  I
also wish to thank Dr Jon Links of the same institute for discussions
and continuing helpful advice.

\pagebreak


\appendix

\section*{Appendix}

Below, we list the nonzero components of our graded R matrices.  The
data are presented in terms of elementary rank $4$ tensors
$e^{ik}_{jl}$, obtained by inserting a copy of $e^k_l$ at each location
of $e^i_j$; where the $e^i_j$ are elementary $2^m \times 2^m$
matrices.  We also use the following notation:

\begin{itemize}
\item
  To increase literacy, we replace $[X]_q$ with $[X]$, we
  often substitute $\overline{q}$ for $q^{-1}$, and we set
  $\Delta \triangleq q - \overline{q}$.

\item
  To convert these \emph{graded} R matrices to the equivalent ungraded
  objects, simply multiply all terms in \textbf{boldface} by $-1$.

\item
  The following notation is a convenient shorthand for the
  frequently appearing \emph{$q$ graded symmetric combination} of rank
  $4$ tensors:
  \begin{eqnarray*}
    q^{x}_{\pm} e^{i k}_{j l}
    & \triangleq &
    q^{x} e^{i k}_{j l} \pm \overline{q}^{x} e^{k i}_{l j}.
  \end{eqnarray*}
  Not using it allows us to present both graded and grading-stripped R
  matrices in one unit.
\end{itemize}

\subsection*{R matrices for $\mathbf{m=1}$}

\input{data/trigrmx1}

\input{data/qrmx1}

\pagebreak

\subsection*{R matrices for $\mathbf{m=2}$}

\input{data/trigrmx2}

\input{data/qrmx2}

\pagebreak

\subsection*{R matrices for $\mathbf{m=3}$}

\input{data/trigrmx3}

\input{data/qrmx3}

\subsection*{R matrices for $\mathbf{m=4}$}

\input{data/trigrmx4}

\input{data/qrmx4}


\bibliographystyle{plain}
\bibliography{DeWit99d}

\end{document}

%% file: data/trigrmx1.tex
Here, $[1]=0$ and $[2]=1$, and $\check{R}^1 (u)$ has $6$ nonzero
components:

\small

\begin{eqnarray*}
  & &
  \hspace{0pt}
  1
  \left\{
    \begin{array}{@{\hspace{0mm}}c@{\hspace{0mm}}}
      e^{1 1}_{1 1}
    \end{array}
  \right\},
  \qquad
  \frac{
    [\alpha + u]
  }{
    [\alpha - u]
  }
  \left\{
    \begin{array}{@{\hspace{0mm}}c@{\hspace{0mm}}}
      e^{2 2}_{2 2}
    \end{array}
  \right\},
  \qquad
  \frac{
    [\alpha]
  }{
    [\alpha - u]
  }
  \left\{
    \begin{array}{@{\hspace{0mm}}c@{\hspace{0mm}}}
      \overline{q}^{u}
      \left\{
      \begin{array}{@{\hspace{0mm}}c@{\hspace{0mm}}}
        e^{1 2}_{1 2}
      \end{array}
      \right\}
      \\
      q^{u}
      \left\{
      \begin{array}{@{\hspace{0mm}}c@{\hspace{0mm}}}
        e^{2 1}_{2 1}
      \end{array}
      \right\}
    \end{array}
  \right\},
  \qquad
  \frac{
    [u]
  }{
    [\alpha - u]
  }
  \left\{
    \begin{array}{@{\hspace{0mm}}c@{\hspace{0mm}}}
      \mathbf{+ 1}
      \left\{
      \begin{array}{@{\hspace{0mm}}c@{\hspace{0mm}}}
        e^{1 2}_{2 1}
      \end{array}
      \right\}
      \\
      - 1
      \left\{
      \begin{array}{@{\hspace{0mm}}c@{\hspace{0mm}}}
        e^{2 1}_{1 2}
      \end{array}
      \right\}
    \end{array}
  \right\}.
\end{eqnarray*}

\normalsize

%% file: data/qrmx1.tex
\noindent
$\check{R}^1$ has $5$ nonzero components:

\small

\begin{eqnarray*}
  \hspace{0pt}
  1
  \left\{
    \begin{array}{@{\hspace{0mm}}c@{\hspace{0mm}}}
       e^{1 1}_{1 1}
    \end{array}
  \right\},
  \qquad
  -
  q^{2 \alpha}
  \left\{
    \begin{array}{@{\hspace{0mm}}c@{\hspace{0mm}}}
      e^{2 2}_{2 2}
    \end{array}
  \right\},
  \qquad
  -
  \Delta
  q^{\alpha}
  [\alpha]
  \left\{
    \begin{array}{@{\hspace{0mm}}c@{\hspace{0mm}}}
      e^{2 1}_{2 1}
    \end{array}
  \right\},
  \qquad
  q^{\alpha}
  \left\{
    \begin{array}{@{\hspace{0mm}}c@{\hspace{0mm}}}
      \mathbf{- 1}
      \left\{
      \begin{array}{@{\hspace{0mm}}c@{\hspace{0mm}}}
        e^{1 2}_{2 1}
      \end{array}
      \right\}
      \\
      + 1
      \left\{
      \begin{array}{@{\hspace{0mm}}c@{\hspace{0mm}}}
        e^{2 1}_{1 2}
      \end{array}
      \right\}
    \end{array}
  \right\}.
\end{eqnarray*}

\normalsize

\vspace{-5mm}

%% file: data/trigrmx2.tex
Here, $[1]=[4]=0$ and $[2]=[3]=1$, and $\check{R}^2 (u)$ has $36$
nonzero components:

\small

\begin{eqnarray*}
  & &
  \hspace{-70pt}
  1
  \left\{
    \begin{array}{@{\hspace{0mm}}c@{\hspace{0mm}}}
      e^{1 1}_{1 1}
    \end{array}
  \right\},
  \qquad
  \frac{[\alpha + u]
  }{
    [\alpha - u]
  }
  \left\{
    \begin{array}{@{\hspace{0mm}}c@{\hspace{0mm}}}
      e^{2 2}_{2 2},
      e^{3 3}_{3 3}
    \end{array}
  \right\},
  \qquad
  \frac{
    [\alpha + u]
    [\alpha + 1 + u]
  }{
    [\alpha - u]
    [\alpha + 1 - u]
  }
  \left\{
    \begin{array}{@{\hspace{0mm}}c@{\hspace{0mm}}}
      e^{4 4}_{4 4}
    \end{array}
  \right\},
  \\[1mm]
  \hline
  \\[1mm]
  & &
  \hspace{-70pt}
  \frac{
    [\alpha]
  }{
    [\alpha - u]
  }
  \left\{
    \begin{array}{@{\hspace{0mm}}c@{\hspace{0mm}}}
      \overline{q}^{u}
      \left\{
      \begin{array}{@{\hspace{0mm}}c@{\hspace{0mm}}}
        e^{1 2}_{1 2},
        e^{1 3}_{1 3}
      \end{array}
      \right\}
      \\
      q^{u}
      \left\{
      \begin{array}{@{\hspace{0mm}}c@{\hspace{0mm}}}
        e^{2 1}_{2 1},
        e^{3 1}_{3 1}
      \end{array}
      \right\}
    \end{array}
  \right\},
  \qquad
  \frac{
    [\alpha + 1]
    [\alpha + u]
  }{
    [\alpha - u]
    [\alpha + 1 - u]
  }
  \left\{
    \begin{array}{@{\hspace{0mm}}c@{\hspace{0mm}}}
      q^{u}
      \left\{
      \begin{array}{@{\hspace{0mm}}c@{\hspace{0mm}}}
        e^{4 3}_{4 3},
        e^{4 2}_{4 2}
      \end{array}
      \right\}
      \\
      \overline{q}^{u}
      \left\{
      \begin{array}{@{\hspace{0mm}}c@{\hspace{0mm}}}
        e^{3 4}_{3 4},
        e^{2 4}_{2 4}
      \end{array}
      \right\}
    \end{array}
  \right\},
  \\
  & &
  \hspace{-70pt}
  \frac{
    [\alpha]
    [\alpha + 1]
  }{
    [\alpha - u]
    [\alpha + 1 - u]
  }
  \left\{
    \begin{array}{@{\hspace{0mm}}c@{\hspace{0mm}}}
      \overline{q}^{ 2 u}
      \left\{
      \begin{array}{@{\hspace{0mm}}c@{\hspace{0mm}}}
        e^{1 4}_{1 4},
      \end{array}
      \right\}
      \\
      q^{2 u}
      \left\{
      \begin{array}{@{\hspace{0mm}}c@{\hspace{0mm}}}
        e^{4 1}_{4 1}
      \end{array}
      \right\}
    \end{array}
  \right\},
  \quad
  \frac{
    1
  }{
    \Delta^{2}
    [\alpha - u]
    [1 + \alpha - u]
  }
  \left\{
    \begin{array}{@{\hspace{0mm}}c@{\hspace{0mm}}}
      f(\overline{q})
      \left\{
      \begin{array}{@{\hspace{0mm}}c@{\hspace{0mm}}}
        e^{2 3}_{2 3}
      \end{array}
      \right\}
      \\
      f(q)
      \left\{
      \begin{array}{@{\hspace{0mm}}c@{\hspace{0mm}}}
        e^{3 2}_{3 2}
      \end{array}
      \right\}
    \end{array}
  \right\},
  \\[1mm]
  \hline
  \\[1mm]
  & &
  \hspace{-70pt}
  \frac{
    [u]
  }{
    [\alpha - u]
  }
  \left\{
    \begin{array}{@{\hspace{0mm}}c@{\hspace{0mm}}}
      \mathbf{+ 1}
      \left\{
      \begin{array}{@{\hspace{0mm}}c@{\hspace{0mm}}}
        e^{1 2}_{2 1},
        e^{1 3}_{3 1}
      \end{array}
      \right\}
      \\
      - 1
      \left\{
      \begin{array}{@{\hspace{0mm}}c@{\hspace{0mm}}}
        e^{2 1}_{1 2},
        e^{3 1}_{1 3}
      \end{array}
      \right\}
    \end{array}
  \right\},
  \qquad
  \frac{
    [u]
    [\alpha + u]
  }{
    [\alpha - u]
    [\alpha + 1 - u]
  }
  \left\{
    \begin{array}{@{\hspace{0mm}}c@{\hspace{0mm}}}
      \mathbf{- 1}
      \left\{
      \begin{array}{@{\hspace{0mm}}c@{\hspace{0mm}}}
        e^{4 3}_{3 4},
        e^{4 2}_{2 4}
      \end{array}
      \right\}
      \\
      + 1
      \left\{
      \begin{array}{@{\hspace{0mm}}c@{\hspace{0mm}}}
        e^{3 4}_{4 3},
        e^{2 4}_{4 2}
      \end{array}
      \right\}
    \end{array}
  \right\},
  \\
  & &
  \hspace{-70pt}
  \frac{
    [u - 1]
    [u]
  }{
    [\alpha - u]
    [\alpha + 1 - u]
  }
  \left\{
    \begin{array}{@{\hspace{0mm}}c@{\hspace{0mm}}}
      e^{1 4}_{4 1}
      \\
      e^{4 1}_{1 4}
    \end{array}
  \right\},
  \qquad
  -
  \frac{
  {[u]}^2
  }{
    [\alpha - u]
    [\alpha + 1 - u]
  }
  \left\{
    \begin{array}{@{\hspace{0mm}}c@{\hspace{0mm}}}
      e^{2 3}_{3 2}
      \\
      e^{3 2}_{2 3}
    \end{array}
  \right\},
  \\[1mm]
  \hline
  \\[1mm]
  & &
  \hspace{-70pt}
  \frac{
    [\alpha]^{\frac{1}{2}}
    [\alpha + 1]^{\frac{1}{2}}
    [u]
  }{
    [\alpha - u]
    [\alpha + 1 - u]
  }
  \left\{
    \begin{array}{@{\hspace{0mm}}c@{\hspace{0mm}}}
      q^{u}
      \left\{
      \begin{array}{@{\hspace{0mm}}c@{\hspace{0mm}}}
        q^{\frac{1}{2}}
        \left\{
        \begin{array}{@{\hspace{0mm}}c@{\hspace{0mm}}}
          \mathbf{- e^{4 1}_{3 2}}
          \\
                  + e^{3 2}_{4 1}
        \end{array}
        \right\},
        \overline{q}^{\frac{1}{2}}
        \left\{
        \begin{array}{@{\hspace{0mm}}c@{\hspace{0mm}}}
          \mathbf{+ e^{4 1}_{2 3}}
          \\
                  - e^{2 3}_{4 1}
        \end{array}
        \right\}
      \end{array}
      \right\}
      \\
      \overline{q}^{u}
      \left\{
      \begin{array}{@{\hspace{0mm}}c@{\hspace{0mm}}}
        \overline{q}^{\frac{1}{2}}
        \left\{
        \begin{array}{@{\hspace{0mm}}c@{\hspace{0mm}}}
          \mathbf{+ e^{1 4}_{2 3}}
          \\
                  - e^{2 3}_{1 4}
        \end{array}
        \right\},
        q^{\frac{1}{2}}
        \left\{
        \begin{array}{@{\hspace{0mm}}c@{\hspace{0mm}}}
          \mathbf{- e^{1 4}_{3 2}}
          \\
                  + e^{3 2}_{1 4}
        \end{array}
        \right\}
      \end{array}
      \right\}
    \end{array}
  \right\},
\end{eqnarray*}
\normalsize
where
$
  f(q)
  =
  - 2 q
  + q^{2 u} (q - \overline{q})
  + q^{2 \alpha} (q + \overline{q})
$.

\vspace{\baselineskip}

%% file: data/qrmx2.tex
\noindent
$\check{R}^2$ has $26$ nonzero components:

\small

\begin{eqnarray*}
  & &
  \hspace{-20pt}
  1
  \left\{
    \begin{array}{@{\hspace{0mm}}c@{\hspace{0mm}}}
       e^{1 1}_{1 1}
    \end{array}
  \right\},
  \qquad
  -
  q^{2 \alpha}
  \left\{
    \begin{array}{@{\hspace{0mm}}c@{\hspace{0mm}}}
      e^{2 2}_{2 2},
      e^{3 3}_{3 3}
    \end{array}
  \right\},
  \qquad
  q^{4 \alpha + 2}
  \left\{
    \begin{array}{@{\hspace{0mm}}c@{\hspace{0mm}}}
       e^{4 4}_{4 4}
    \end{array}
  \right\},
  \\[1mm]
  \hline
  \\[1mm]
  & &
  \hspace{-20pt}
  -
  \Delta
  q^{\alpha}
  [\alpha]
  \left\{
  \begin{array}{@{\hspace{0mm}}c@{\hspace{0mm}}}
    e^{2 1}_{2 1},
    e^{3 1}_{3 1}
  \end{array}
  \right\},
  \quad
  \Delta
  q^{3 \alpha + 1}
  [\alpha + 1]
  \left\{
  \begin{array}{@{\hspace{0mm}}c@{\hspace{0mm}}}
    e^{4 3}_{4 3},
    e^{4 2}_{4 2}
  \end{array}
  \right\},
  \quad
  \Delta^2
  q^{2 \alpha + 1}
  [\alpha]
  [\alpha + 1]
  \left\{
  \begin{array}{@{\hspace{0mm}}c@{\hspace{0mm}}}
    e^{4 1}_{4 1}
  \end{array}
  \right\},
  \quad
  \Delta
  q^{2 \alpha + 1}
  \left\{
  \begin{array}{@{\hspace{0mm}}c@{\hspace{0mm}}}
    e^{3 2}_{3 2}
  \end{array}
  \right\},
  \\[1mm]
  \hline
  \\[1mm]
  & &
  \hspace{-20pt}
  q^{\alpha}
  \left\{
  \begin{array}{@{\hspace{0mm}}c@{\hspace{0mm}}}
    \mathbf{- 1}
    \left\{
    \begin{array}{@{\hspace{0mm}}c@{\hspace{0mm}}}
      e^{1 2}_{2 1},
      e^{1 3}_{3 1}
    \end{array}
    \right\}
    \\
    + 1
    \left\{
    \begin{array}{@{\hspace{0mm}}c@{\hspace{0mm}}}
      e^{2 1}_{1 2},
      e^{3 1}_{1 3}
    \end{array}
    \right\}
  \end{array}
  \right\},
  \qquad
  q^{3 \alpha + 1}
  \left\{
  \begin{array}{@{\hspace{0mm}}c@{\hspace{0mm}}}
    \mathbf{- 1}
    \left\{
    \begin{array}{@{\hspace{0mm}}c@{\hspace{0mm}}}
      e^{4 3}_{3 4},
      e^{4 2}_{2 4}
    \end{array}
    \right\}
    \\
    + 1
    \left\{
    \begin{array}{@{\hspace{0mm}}c@{\hspace{0mm}}}
      e^{3 4}_{4 3},
      e^{2 4}_{4 2}
    \end{array}
    \right\}
  \end{array}
  \right\},
  \qquad
  q^{2 \alpha}
  \left\{
  \begin{array}{@{\hspace{0mm}}c@{\hspace{0mm}}}
    e^{1 4}_{4 1}
    \\
    e^{4 1}_{1 4}
  \end{array}
  \right\},
  \qquad
  -
  q^{2 \alpha + 1}
  \left\{
  \begin{array}{@{\hspace{0mm}}c@{\hspace{0mm}}}
    e^{2 3}_{3 2}
    \\
    e^{3 2}_{2 3}
  \end{array}
  \right\},
  \\[1mm]
  \hline
  \\[1mm]
  & &
  \hspace{-20pt}
  \Delta
  q^{2 \alpha + 1}
  [\alpha]^{\frac{1}{2}}
  [\alpha + 1]^{\frac{1}{2}}
  \left\{
    q^{\frac{1}{2}}
    \left\{
    \begin{array}{@{\hspace{0mm}}c@{\hspace{0mm}}}
      \mathbf{- e^{4 1}_{3 2}}
      \\
              + e^{3 2}_{4 1}
    \end{array}
    \right\},
    \overline{q}^{\frac{1}{2}}
    \left\{
    \begin{array}{@{\hspace{0mm}}c@{\hspace{0mm}}}
      \mathbf{+ e^{4 1}_{2 3}}
      \\
              - e^{2 3}_{4 1}
    \end{array}
    \right\}
  \right\}.
\end{eqnarray*}

\normalsize

%% file: data/trigrmx3.tex
\noindent
Here, $[i]=0$ for $i\in \{1;5,6,7\}$ and $[i]=1$ for $i\in \{2,3,4;8\}$.
The reader will have by now appreciated the recurring patterns
in the components of our R matrices. To save space, we introduce
a little more notation, which eliminates the $q$ brackets altogether:
\begin{eqnarray*}
  S^\pm_i
  & \triangleq &
  [ \alpha + i \pm u ]_q,
  \\
  A^z_i
  & \triangleq &
  [ \alpha + i ]^z_q,
  \qquad
  \mathrm{where~} z\in \{{\textstyle \frac{1}{2}}, 1\},
  \\
  U^z_i
  & \triangleq &
  [ u - i ]^z_q,
  \qquad
  \mathrm{where~} z\in \{1, 2\},
\end{eqnarray*}
and $i \in \{ 0,1, \dots, m-1\}$.
With this notation,
$\check{R}^3 (u)$ has $216$ nonzero components:

\small

\begin{eqnarray*}
  & &
  \hspace{-20pt}
  1
  \left\{

      \right\}
    \end{array}
    \right\}
    \end{array}
  \right\},
\end{eqnarray*}
\normalsize
where:
\small
\begin{eqnarray*}
  f_1 (q)
  & = &
  - 2 \overline{q}
  + (q^{1 + 2 \alpha} + \overline{q}^{1 + 2 \alpha})
  - \overline{q}^{2 u} (q - \overline{q})
  \\
  f_2 (q)
  & = &
  - 2 q
  + (\overline{q}^{3 + 2 \alpha} + q^{3 + 2 \alpha})
  + q^{2 u} (q - \overline{q})
  \\
  f_3 (q)
  & = &
  - \overline{q}^{u}
  (
    2 \overline{q}^{2}
    - (q^{2 + 2 \alpha} + \overline{q}^{2 + 2 \alpha})
    + \overline{q}^{2 u} (q^{2} - \overline{q}^{2})
  )
  \\
  f_4 (q)
  & = &
  \overline{q}^{u}
  (
    - 2 
    + (q^{2 + 2 \alpha} + \overline{q}^{2 + 2 \alpha})
    - (q^{2 - 2 u} + \overline{q}^{2 - 2 u})
    + (q^{2u} + \overline{q}^{2u})
  )
  \\
  f_5 (q)
  & = &
  \overline{q}^{u}
  (
    - 2 q^{2}
    + (q^{2 + 2 \alpha} + \overline{q}^{2 + 2 \alpha})
    + q^{2 u} (q^{2} - \overline{q}^{2})
  )
  \\
  f_6 (q)
  & = &
  q (q + \overline{q})
  - (q^{2 + 2 \alpha} + \overline{q}^{2 + 2 \alpha})
  - q^{2 u - 1} (q - \overline{q}).
\end{eqnarray*}

\normalsize

%% file: data/qrmx3.tex
\noindent
$\check{R}^3$ has $139$ nonzero components:

\small

\begin{eqnarray*}
  & &
  \hspace{-20pt}
  1
  \left\{

    \right\}
  \right\}.
\end{eqnarray*}

\normalsize

%% file: data/trigrmx4.tex
\noindent
Here, $[i]$ is $0$ for $i \in \{1;6,7,8,9,10,11;16\}$, and $1$ for $i
\in \{2,3,4,5;12,13,14,15\}$, $\check{R}^4 (u)$ has $1296$ nonzero
components:

\tiny

\begin{eqnarray*}
  &&
  \hspace{-20pt}
  1
  \left\{
    %
        \right\}
      \end{array}
      \right\}
    \end{array}
  \right\},
\end{eqnarray*}
\normalsize
where:
\tiny
\begin{eqnarray*}
  \hspace{0pt}
  g_{ 1} (q)
  & = &
  - 2 q
  + (q^{1 + 2 \alpha} + \overline{q}^{1 + 2 \alpha})
  + q^{2 u} (q - \overline{q})
  \\
  \hspace{0pt}
  g_{ 2} (q)
  & = &
  \overline{q}^{u}
  (
    - 2
    + (q^{2 + 2 \alpha} + \overline{q}^{2 + 2 \alpha})
    + (q^{2 u} + \overline{q}^{2 u})
    - (q^{2 u - 2} + \overline{q}^{2 u - 2})
  )
  \\
  \hspace{0pt}
  g_{ 3} (q)
  & = &
  - 2 q^{2}
  + (q^{2\alpha+2} + \overline{q}^{2\alpha+2})
  + q^{2 u} (q^{2} - \overline{q}^2)
  \\
  \hspace{0pt}
  g_{ 4} (q)
  & = &
  - 2 q
  + (q^{3 + 2 \alpha} + \overline{q}^{3 + 2 \alpha})
  + q^{2 u} (q - \overline{q})
  \\
  \hspace{0pt}
  g_{ 5} (q)
  & = &
  - 2 q^{2}
  + (q^{4 + 2 \alpha} + \overline{q}^{4 + 2 \alpha})
  + q^{2 u} (q^{2} - \overline{q}^{2})
  \\
  \hspace{0pt}
  g_{ 6} (q)
  & = &
  - 2 q^{3}
  + (q^{3 + 2 \alpha} + \overline{q}^{3 + 2 \alpha})
  + q^{2 u} (q^{3} - \overline{q}^{3})
  \\
  \hspace{0pt}
  g_{ 7} (q)
  & = &
  - 2 q
  + (q^{3 + 2 \alpha} + \overline{q}^{3 + 2 \alpha})
  + q (q^{2 u} + \overline{q}^{2 u})
  - (q^{2 u - 3} + \overline{q}^{2 u - 3})
  \\
  \hspace{0pt}
  g_{ 8} (q)
  & = &
  q^{u}
  (
    - 2
    +  (q^{4 + 2 \alpha} + \overline{q}^{4 + 2 \alpha})
    +  (q^{2 u} + \overline{q}^{2 u})
    -  ( q^{2 u - 2} + \overline{q}^{2 u - 2})
  )
  \\
  \hspace{0pt}
  g_{ 9} (q)
  & = &
  - 2 q
  + (q^{5 + 2 \alpha} - \overline{q}^{5 + 2 \alpha})
  + q^{2 u} (q - \overline{q})
  \\
  \hspace{0pt}
  g_{10} (q)
  & = &
  1
  + 2 q^{2}
  + 3 q^{4}
  - 2 \overline{q}^{2 + 2 \alpha}
  - 2 \overline{q}^{2 \alpha}
  - 2 q^{4 + 2 \alpha}
  - 2 q^{6 + 2 \alpha}
  + (q^{6 + 4 \alpha} + \overline{q}^{6 + 4 \alpha})
  + 2 q^{2 u}
  -   q^{4 u}
  - q^{2 u} (q^{2} - \overline{q}^{2})
  \\
  \hspace{0pt}
  & &
  - 2 q^{4 + 2 u}
  + q^{4 u - 4}
  + q^{2 u} (q^{4 + 2 \alpha} - \overline{q}^{4 + 2 \alpha})
  + q^{4 u} (q^{2} - \overline{q}^{2})
  + q^{2 u} (q^{6 + 2 \alpha} - \overline{q}^{6 + 2 \alpha})
  - q^{2 u} (q^{2 \alpha} - \overline{q}^{2 \alpha})
  \\
  \hspace{0pt}
  & &
  - q^{2 u} (q^{2 + 2 \alpha} - \overline{q}^{2 + 2 \alpha})
  \\
  \hspace{0pt}
  g_{11} (q)
  & = &
    3
  -   \overline{q}^2
  + 5 q^{2}
  -   q^{4}
  +   (q^{6 + 4 \alpha} + \overline{q}^{6 + 4 \alpha})
  - 4 \overline{q}^{2 + 2 \alpha}
  - 4 q^{4 + 2 \alpha}
  + 2 q^{- 2 - 2 \alpha + 2 u}
  + 2 q^{4 + 2 \alpha + 2 u}
  \\
  \hspace{0pt}
  & &
  -   (q^{2 \alpha + 2 u} + \overline{q}^{2 \alpha + 2 u})
  -   (q^{2 + 2 \alpha + 2 u} - \overline{q}^{2 + 2 \alpha + 2 u})
  +   (q^{4 + 2 \alpha - 2 u} - \overline{q}^{4 + 2 \alpha - 2 u})
  -   (q^{6 + 2 \alpha - 2 u} + \overline{q}^{6 + 2 \alpha - 2 u})
  \\
  \hspace{0pt}
  & &
  - 2 q^{2} (q^{2 u} + \overline{q}^{2 u})
  + 2 q^{4 - 2 u}
  - 2 \overline{q}^{2 - 4 u}
  - 2 q^{2 u}
  +   q^{4 u}
  +   \overline{q}^{4 - 4 u}
  + 4 \overline{q}^{2 - 2 u}
  \\
  \hspace{0pt}
  g_{12} (q)
  & = &
  6
  + (q^{2} + \overline{q}^{2})
  - (q^{4} + \overline{q}^{4})
  - 2 (q^{2 + 2 \alpha} + \overline{q}^{2 + 2 \alpha})
  - 2 (q^{4 + 2 \alpha} + \overline{q}^{4 + 2 \alpha})
  +   (q^{6 + 4 \alpha} + \overline{q}^{6 + 4 \alpha})
  + 2 q^{4 - 2 u},
\end{eqnarray*}
\normalsize
and:
\tiny
\begin{eqnarray*}
  \hspace{-90pt}
  h_{ 1} (q)
  & = &
  - q (q + \overline{q})
  + (q^{2 + 2 \alpha} + \overline{q}^{2 + 2 \alpha})
  + q^{2 u - 1} (q - \overline{q})
  \\
  \hspace{-90pt}
  h_{ 2} (q)
  & = &
  - (q^{3 + 2 \alpha} + \overline{q}^{3 + 2 \alpha})
  - q^{u + 1} (q^{u} - \overline{q}^{u})
  + q^{u} (q^{3 + u} + \overline{q}^{3 + u})
  \\
  \hspace{-90pt}
  h_{ 3} (q)
  & = &
  (q + \overline{q})
  - (q^{3 + 2 \alpha} + \overline{q}^{3 + 2 \alpha})
  - (q^{1 - 2 u} + \overline{q}^{1 - 2 u})
  + (q^{3 - 2 u} + \overline{q}^{3 - 2 u})
  \\
  \hspace{-90pt}
  h_{ 4} (q)
  & = &
  - 2 q^{2}
  + q (
        \overline{q}^{3 + 2 \alpha}
        +
        q^{3 + 2 \alpha}
      )
  + q^{2 u - 1} (q - \overline{q})
  \\
  \hspace{-90pt}
  h_{ 5} (q)
  & = &
    q^{2 u - 1} (q - \overline{q})
  - q^{2} (q^{2} - \overline{q}^{2})
  + q (q^{3 + 2 \alpha} + \overline{q}^{3 + 2 \alpha})
  \\
  \hspace{-90pt}
  h_{ 6} (q)
  & = &
  - (q^{4 + 2 \alpha} + \overline{q}^{4 + 2 \alpha})
  - q^{u} (q^{u} - \overline{q}^{u})
  + q^{u} (q^{2 - u} + \overline{q}^{2 - u}).
\end{eqnarray*}

\normalsize

%% file: data/qrmx4.tex
\noindent
$\check{R}^4$ has $758$ nonzero components:

\tiny

\begin{eqnarray*}
  &&
  \hspace{-20pt}
  1
  \left\{

      \right\}
    \end{array}
  \right\}.
\end{eqnarray*}

\normalsize